\documentclass[12pt]{iopart}
\usepackage{graphicx}
\begin{document}

\title[A convergent scheme for initial shape identification]{An unconditionally convergent iterative scheme for initial shape identification in small deformations}

\author{M. Sellier}

\address{Schott AG, Hattenbergstrasse 10, 55122 Mainz, Germany}

\ead{mathieu.sellier@itwm.fraunhofer.de}

\begin{abstract}
The question of interest in the present study is: ``Given a body subject to mechanical loads, how to define the initial geometry so that the deformed one matches precisely a prescribed shape?'' This question is particularly relevant in forming processes where the tolerated mismatch between the deformed and desired geometries may be lower than a Micron. The method proposed here uses as a first ``guess'' to the required initial geometry the desired one, then it updates iteratively the locations of a set of boundary points so that their locations in the deformed configuration come closer and closer to the desired ones. The scheme is shown to converge unconditionally for small deformations in the sense that arbitrarily small mismatch between the deformed and desired shape can be achieved. Moreover, since it is based entirely on geometric considerations, the convergence should not be affected by the nature of material, i.e. it is independent of the constitutive law. The success of the method is illustrated by considering an example.     
\end{abstract}

\maketitle

\section{Introduction}
The manufacturing of high-precision glass or metal products is a challenge due to the complexity of the phenomena involved during the forming process and the difficulty to control in practice all the relevant parameters. A question engineers often have to face is: ``given a set of loads applied to a solid body, what should be the initial geometry so that the end geometry (after deformation) matches the prescribed one?'' The development of efficient and accurate numerical techniques such as the Finite Element method allowed significant progress towards the solution to this inverse problem. A popular approach consists in parameterizing somehow the contour of the work-piece, defining an objective function related to the mismatch between the deformed and desired work-piece geometries and applying the extensive theory of optimization in order to reduce the objective function. This approach has proven its applicability for the problem of tool design in the metal forging process (see \cite{Chung03,Vielledent01,Sousa02} and reference therein). Many authors favour ``gradient-based'' optimization techniques and the condition to the success of these is the accurate and efficient evaluation of the sensitivities. The sensitivity analysis may be performed using the Finite Difference method, the direct differentiation method, \cite{Vielledent01,Sousa02}, or the adjoint-state method, \cite{Chung03,Sprekels}. In many practical cases, the pragmatic user of a commercial code for the computation of the deformations is reduced to the use of  the Finite Difference method to compute the sensitivities and this approach is known for its lack of accuracy and efficiency, \cite{Chung03}.

An alternative approach which overcome this issue was first introduced by Park \etal, \cite{Park83}. It consists in starting from the desired work-piece geometry and back-tracking the deformation path. The idea was applied with success by various authors to preform or die design in metal forging, \cite{Kim90,Zhao94}. A similar concept was exploited in \cite{Sellier04} in combination with an iterative scheme to update the locations of a set boundary points in order to identify the required initial geometry. This method proved its applicability and its efficiency in the context of high-precision glass forming and was easily implemented in a commercial code. In spite of the non-linear viscoelastic behaviour of the glass, only few iterations were required to reduce the mismatch between the deformed and desired work-piece geometries below a Micron. A formal analysis explained the good convergence properties of this iterative scheme and revealed a weakness. It is, in theory, only applicable to deformations dominated by shear or volumetric strains.

The aim of the present work is to present a modified scheme which circumvent this restriction. The main idea remains essentially the same but a corrective term is added for the update of the boundary points locations. The next two sections are devoted to the description of the scheme and the analysis of its convergence properties. The scheme is next applied to a case where the prescribed displacement field combines shear and volumetric strains in equal magnitude in order to assess its performance.  
\section{Description of the iterative scheme}
\begin{figure}
\centering
\scalebox{0.5}{\includegraphics{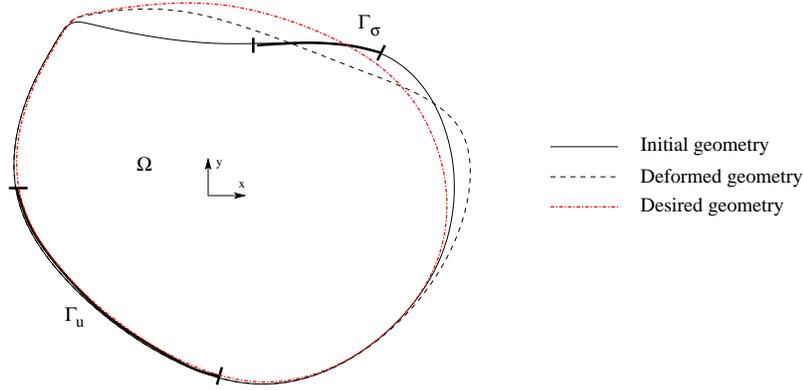}}
\caption{Sketch of the general forward mechanical problem.}
\label{fig:geom}
\end{figure}

We restrict the analysis to a mechanical problem in a two-dimensional domain $\Omega$ such that the displacement field $\vec{U}$ is a two-dimensional vector whose components depend only on the coordinates $(x,y)$ of the Cartesian coordinate system associated with the body shown on Figure \ref{fig:geom}. 
The body is subject to mechanical loads. On $\Gamma _u$, displacements are prescribed and traction is imposed on $\Gamma _{\sigma}$. Clearly, the forward problem corresponds to calculating the deformed shape given the initial one and the boundary conditions while the inverse problem we wish to solve is: ``given a desired shape, find the initial one such that the deformed work-piece geometry matches the desired one''. We will assume in the following that the solution to the forward problem is obtainable by analytical or, more likely, numerical means. Consequently, the exact or approximated displacement field can be found for all points in $\Omega$.

\begin{figure}
\centering
\scalebox{0.4}{\includegraphics{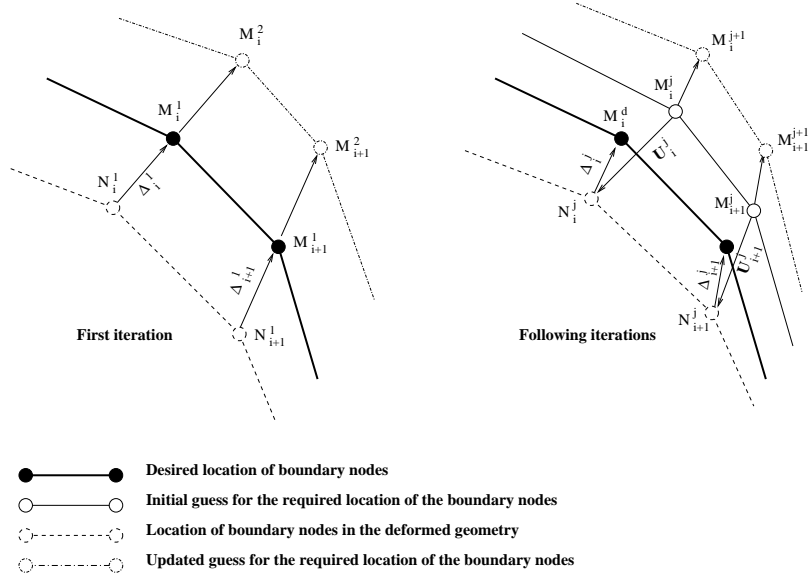}}
\caption{The iterative algorithm with the corrective term $\vec{B}_i^j$ set to zero.}
\label{sellierfig:2}
\end{figure}

In order to describe the algorithm which tackles the inverse problem of identifying the required initial geometry, notations are first detailed. A number of points are introduced at the boundary. If a discrete method is used to solve the forward problem, these points naturally correspond to the boundary nodes of the mesh. Let $M_1^d,\ldots,M_L^d$ denote the $L$ boundary nodes of the desired work-piece geometry and $M_1^{ini},\ldots,M_L^{ini}$ denote the associated $L$ boundary nodes of the required initial geometry (solution to the inverse problem). The $M_i^{ini}$ are found iteratively as the limit of a series of $M_i^j$, where the superscript $j$ denotes the iteration number. At each iteration, $N_i^j$ corresponds to the new location of the node $M_i^j$ in the deformed geometry and $\vec{U}_i^j$ is the associated displacement. The algorithm is best described in pseudo-code notation as follows:
\begin{enumerate}
\item \underline{First iteration}: use $M_i^d\;,\;i\in[1,L]$ as a first guess to $M_i^{ini}\;,\;i\in[1,L]$ and calculate the \textsl{residual vector} $\vec{\Delta} _i^1\;,\;i\in[1,L]$ whose norm gives a measure of how far the node in the deformed geometry is from the desired location:\\
\texttt{for i=1 to L} \{\\
\begin{eqnarray}
\label{eq:sub1}
\vec{OM}_i^1 & = & \vec{OM}_i^d \; ; \\
\label{eq:sub2}
\vec{ON}_i^1 & = & \vec{OM}_i^1+\vec{U}_i^1 \; ; \\
\label{eq:sub3}
\vec{\Delta}_i^1 & = & \vec{OM}_i^d-\vec{ON}_i^1 \; ;
\end{eqnarray}
\indent \indent \} \texttt{j=2};
\item \underline{Following iterations}: update the previous estimate to $M_i^{ini}\;,\;i\in[1,L]$ and calculate the residual vector $\vec{\Delta} _i^j\;,\;i\in[1,L]$:\\
\texttt{Do} \{\\
\indent \indent \texttt{for i=1 to L} \{\\
\begin{eqnarray}
\label{eq:sub4}
\vec{OM}_i^j & = & \vec{OM}_i^{j-1}+\vec{\Delta}_i^{j-1}+\vec{B}_i^j \; ; \\
\label{eq:sub5}
\vec{ON}_i^j & = & \vec{OM}_i^j+\vec{U}_i^j \; ; \\
\label{eq:sub6}
\vec{\Delta}_i^j & = & \vec{OM}_i^d-\vec{ON}_i^j \; ;
\end{eqnarray}
\indent \indent \} \texttt{j=j+1};\}\\
\indent \indent \texttt{While $\max(||\vec{\Delta} _i^j||)>\epsilon$}
\end{enumerate}
Restated in simple terms, the initial guess for the required initial boundary node locations is chosen to be the location of the nodes of the desired geometry. At each iteration the residual vector ($\vec{\Delta}_i^j$) whose norm measures how far the deformed geometry is from the desired one is evaluated and added to the previous guess of the required initial boundary node location. $\vec{B}_i^j$ is a corrective term which enlarges the range of applicability of the scheme and shall be defined subsequently. If this corrective term vanishes, the intuitive scheme proposed in \cite{Sellier04} and illustrated on Figure \ref{sellierfig:2} is recovered. Since this scheme is entirely based on geometric consideration, it does not depend on the material behaviour and should therefore be applicable regardless of the constitutive law. This was shown in \cite{Sellier04} where the scheme was used for the identification of the required initial geometry of a glass piece, a strongly nonlinear material with stress and structure relaxation.
\section{Convergence analysis}
In order to assess the convergence properties of the scheme (i.e. will $N_i^j\;,\;i\in[1,L]$ converge to $M_i^d\;,\;i\in[1,L]$ for increasing $j$, and if so, how fast?), we perform the same analysis as that presented in \cite{Sellier04}. Firstly, it is necessary to note that as a consequence of the algorithm,
\begin{equation}
\vec{\Delta}_i^j = \vec{U}_i^{j-1}-\vec{U}_i^j-\vec{B}_i^j \; .
\label{eq:resvec}
\end{equation} 
This result is proven by simple algebraic manipulations as follows;\\
Using eq.\ (\ref{eq:sub6}) gives 
$$\vec{\Delta}_i^j-\vec{\Delta}_i^{j-1}=\vec{OM}_i^d-\vec{ON}_i^j-\left( \vec{OM}_i^d-\vec{ON}_i^{j-1} \right) \; .$$
Substitution of eq.\ (\ref{eq:sub4}) into eq.\ (\ref{eq:sub5}) and of eq.\ (\ref{eq:sub5}) into the previous equation yields
$$\vec{\Delta}_i^j-\vec{\Delta}_i^{j-1}=\vec{OM}_i^{j-2}+\vec{\Delta} _i^{j-2}+\vec{B}_i^{j-1}+\vec{U}_i^{j-1}-\left( \vec{OM}_i^{j-1}+\vec{\Delta} _i^{j-1}+\vec{B}_i^j+\vec{U}_i^j \right) \; ,$$
but since according to eq.\ (\ref{eq:sub4}), $\vec{OM}_i^{j-1}=\vec{OM}_i^{j-2}+\vec{\Delta}_i^{j-2}+\vec{B}_i^{j-1}$, eq.\ (\ref{eq:resvec}) is recovered which proves the result.

At iteration $(j-1)$, the displacement of the node $M_i^{j-1}$ with coordinates $(x_i^{j-1},y_i^{j-1})$ can be found by Taylor series expansion around $M_i^d$ with coordinates $(x_i^d,y_i^d)$. Accordingly,
\begin{eqnarray}
\label{eq:urijm1} 
{U_x}_i^{j-1}  & = & U_x \left( x_i^d,y_i^d \right) +\frac{\partial U_x}{\partial x}|_i^d \left( x_i^{j-1}-x_i^d \right)+\frac{\partial U_x}{\partial y}|_i^d \left( y_i^{j-1}-y_i^d \right)+h.o.t. \; , \\
\label{eq:uzijm1}
{U_y}_i^{j-1}  & = & U_y \left( x_i^d,y_i^d \right)+\frac{\partial U_y}{\partial x}|_i^d \left( x_i^{j-1}-x_i^d \right)+\frac{\partial U_y}{\partial y}|_i^d \left( y_i^{j-1}-y_i^d \right)+h.o.t. \; ,
\end{eqnarray}
where $\left( {U_x}_i^{j-1},{U_y}_i^{j-1}\right)$ are the components of the displacement vector $\vec{U}_i^{j-1}$ and it is understood that the spatial derivatives of the displacement are taken at $M_i^d$.\\
Similarly, at iteration $(j)$, the displacement of the node $M_i^j$ can be written as follows:
\begin{eqnarray}
\label{eq:urij} 
{U_x}_i^{j} & = & U_x \left( x_i^d,y_i^d \right) +\frac{\partial U_x}{\partial x}|_i^d \left( x_i^{j}-x_i^d \right)+\frac{\partial U_x}{\partial y}|_i^d \left( y_i^{j}-y_i^d \right)+h.o.t. \; , \\
\label{eq:uzij}
{U_y}_i^{j} & = & U_y \left( x_i^d,y_i^d \right)+\frac{\partial U_y}{\partial x}|_i^d \left( x_i^{j}-x_i^d \right)+\frac{\partial U_y}{\partial y}|_i^d \left( y_i^{j}-y_i^d \right)+h.o.t. \; .
\end{eqnarray}
Subtracting eq.\ (\ref{eq:urij}) to eq.\ (\ref{eq:urijm1}) and eq.\ (\ref{eq:uzij}) to eq.\ (\ref{eq:uzijm1}) gives,
\begin{eqnarray}
{U_x}_i^{j-1}-{U_x}_i^{j} & = & \frac{\partial U_x}{\partial x}|_i^d \left( x_i^{j-1}-x_i^{j} \right)+\frac{\partial U_x}{\partial y}|_i^d \left( y_i^{j-1}-y_i^{j} \right) + h.o.t. \; , \\
{U_y}_i^{j-1}-{U_y}_i^{j} & = & \frac{\partial U_y}{\partial x}|_i^d \left( x_i^{j-1}-x_i^{j} \right)+\frac{\partial U_y}{\partial y}|_i^d \left( y_i^{j-1}-y_i^{j} \right) + h.o.t. \; .
\end{eqnarray}
Moreover, taking eqs.\ (\ref{eq:sub4}) and (\ref{eq:resvec}) into account and remembering that $\frac{\partial U_x}{\partial x} = \epsilon _{xx}$ and $\frac{\partial U_y}{\partial y} = \epsilon _{yy}$, these equations may be rewritten as,
\begin{eqnarray}
\label{eq:deltarf}
{\Delta _x}_i^j & = & -{\epsilon _{xx}}_i^d\left({\Delta _x}_i^{j-1}+{B_x}_i^j\right) -\frac{\partial U_x}{\partial y}|_i^d\left({\Delta _y}_i^{j-1}+{B_y}_i^j\right)-{B_x}_i^j + h.o.t. \; , \\
\label{eq:deltazf}
{\Delta _y}_i^j & = & -\frac{\partial U_y}{\partial x}|_i^d\left({\Delta _x}_i^{j-1}+{B_x}_i^j \right) -{\epsilon _{yy}}_i^d\left({\Delta _y}_i^{j-1} +{B_y}_i^j\right)-{B_y}_i^j + h.o.t. \; ,
\end{eqnarray}
where $({\Delta _x}_i^j,{\Delta _y}_i^j)$ and $({B_x}_i^j,{B_y}_i^j)$ are the x-y components of the residual vector $\vec{\Delta}_i^j$ and the corrective term $\vec{B}_i^j$ respectively. If the corrective term vanishes, it is easy to show (see \cite{Sellier04}) that two extreme cases should produce convergent schemes. The first one corresponds to deformations dominated by volumetric strains (i.e. $\frac{\partial U_x}{\partial y} \ll \epsilon _{xx} < 1$ and $\frac{\partial U_y}{\partial x} \ll \epsilon _{yy} <1$) while the second is applicable when the deformations are dominated by the shear terms, i.e $\epsilon _{xx} \ll \frac{\partial U_x}{\partial y}$ and $\epsilon _{yy} \ll \frac{\partial U_y}{\partial x}$. The scheme without corrective term will be referred to as \textbf{scheme (I)} in the following. A suitable choice for the corrective term $\vec{B}_i^j$ should allow an unconditional convergence. One immediately sees that if $\vec{B}_i^j$ is chosen such that,
\begin{eqnarray}
\label{eq:correctx}
-{\epsilon _{xx}}_i^d{B_x}_i^j-\frac{\partial U_x}{\partial y}|_i^d\left( {\Delta _y}_i^{j-1}+{B_y}_i^j \right)-{B_x}_i^j & = & 0 \; , \\
\label{eq:correcty}
-\frac{\partial U_y}{\partial x}|_i^d\left({\Delta _x}_i^{j-1}+{B_x}_i^j \right)-{\epsilon _{yy}}_i^d{B_y}_i^j-{B_y}_i^j & = & 0 \; ,
\end{eqnarray} 
equations (\ref{eq:deltarf}) and (\ref{eq:deltazf}) simply reduce to
\begin{equation}
{\Delta _x}_i^j = -{\epsilon _{xx}}_i^d{\Delta _x}_i^{j-1} \quad \textrm{and} \quad {\Delta _y}_i^j = -{\epsilon _{yy}}_i^d{\Delta _y}_i^{j-1} \; ,
\label{eq:convscheme}
\end{equation}
which clearly defines a convergent scheme providing the absolute value of the strains is smaller than one which is necessarily the case in small deformations. Note however that since the convergence analysis is based on the Taylor expansion of the displacement field around $M_i^d$ and that only the first order terms are retained, its validity is resctricted to situations when $M_i^j$ is sufficiently close to $M_i^d$. Equations (\ref{eq:correctx}) and (\ref{eq:correcty}) are rewritten in a more readable form as follows,
\begin{eqnarray}
\label{eq:correctxnice}
\left( {\epsilon _{xx}}_i^d+1\right){B_x}_i^j+\frac{\partial U_x}{\partial y}|_i^d{B_y}_i^j & = & -\frac{\partial U_x}{\partial y}|_i^d{\Delta _y}_i^{j-1} \; , \\
\label{eq:correctynice}
\frac{\partial U_y}{\partial x}|_i^d{B_x}_i^j+\left( {\epsilon _{yy}}_i^d+1\right){B_y}_i^j & = & -\frac{\partial U_y}{\partial x}|_i^d{\Delta _x}_i^{j-1} \; .
\end{eqnarray}
Since ${\epsilon _{xx}}_i^d$, ${\epsilon _{yy}}_i^d$, $\frac{\partial U_x}{\partial y}|_i^d$ and $\frac{\partial U_y}{\partial x}|_i^d$ must be evaluated at $M_i^d$, these can be computed in practice at the first iteration and therefore the system of equations (\ref{eq:correctxnice}) and (\ref{eq:correctynice}) can be solved for the two components of the corrective vector $\vec{B}_i^j$. We denote by \textbf{scheme (II)} the one with the corrective term which satisfies equations (\ref{eq:correctxnice}) and (\ref{eq:correctynice}) exactly.\\
Finally, a simple analysis of the orders of magnitude reveals that providing ${\epsilon _{xx}}_i^d \ll 1$, ${\epsilon _{yy}}_i^d \ll 1$, $\frac{\partial U_x}{\partial y}|_i^d \ll 1$ and $\frac{\partial U_y}{\partial x}|_i^d \ll 1$, the corrective vector reduces to $\vec{B}_i^j \simeq \left( -\frac{\partial U_x}{\partial y}|_i^d{\Delta _y}_i^{j-1}, -\frac{\partial U_y}{\partial x}|_i^d{\Delta _x}_i^{j-1}\right)$ and the latter scheme will be referred to as \textbf{scheme (III)}.
\section{Test problem}
As a test case, a known and prescribed displacement field is considered and we seek the initial geometry so that the deformed one matches a disc of radius 0.01. The following displacement field is imposed:
\begin{equation}
\vec{U} = \left( U_x,U_y \right) = \left( \alpha(x+y),\alpha(x-y)\right) \; .
\label{eq:displacement}
\end{equation}
Of course, in this case an analytical solution to the inverse problem can be found analytically by solving the following system of equations,
\begin{eqnarray}
\label{eq:truex}
x_i^d-x_i^{ini} & = & \alpha \left( x_i^{ini}+y_i^{ini}\right) \; , \\
\label{eq:truey}
y_i^d-y_i^{ini} & = & \alpha \left( x_i^{ini}-y_i^{ini}\right) \; , 
\end{eqnarray}
for $(x_i^{ini},y_i^{ini})$ with $(x_i^d,y_i^d)$ belonging to the circle of radius 0.01. This displacement field is however a good candidate to assess the proposed method since the terms $\epsilon _{xx}$, $\epsilon _{yy}$, $\frac{\partial U_x}{\partial y}$ and $\frac{\partial U_y}{\partial x}$ will all have the same magnitude equal to $\alpha$ and this is precisely the situation when, according to the previous analysis, the convergence of \textbf{scheme (I)} can not be guaranteed. 
\begin{figure}
\begin{center}
\scalebox{0.7}{\includegraphics{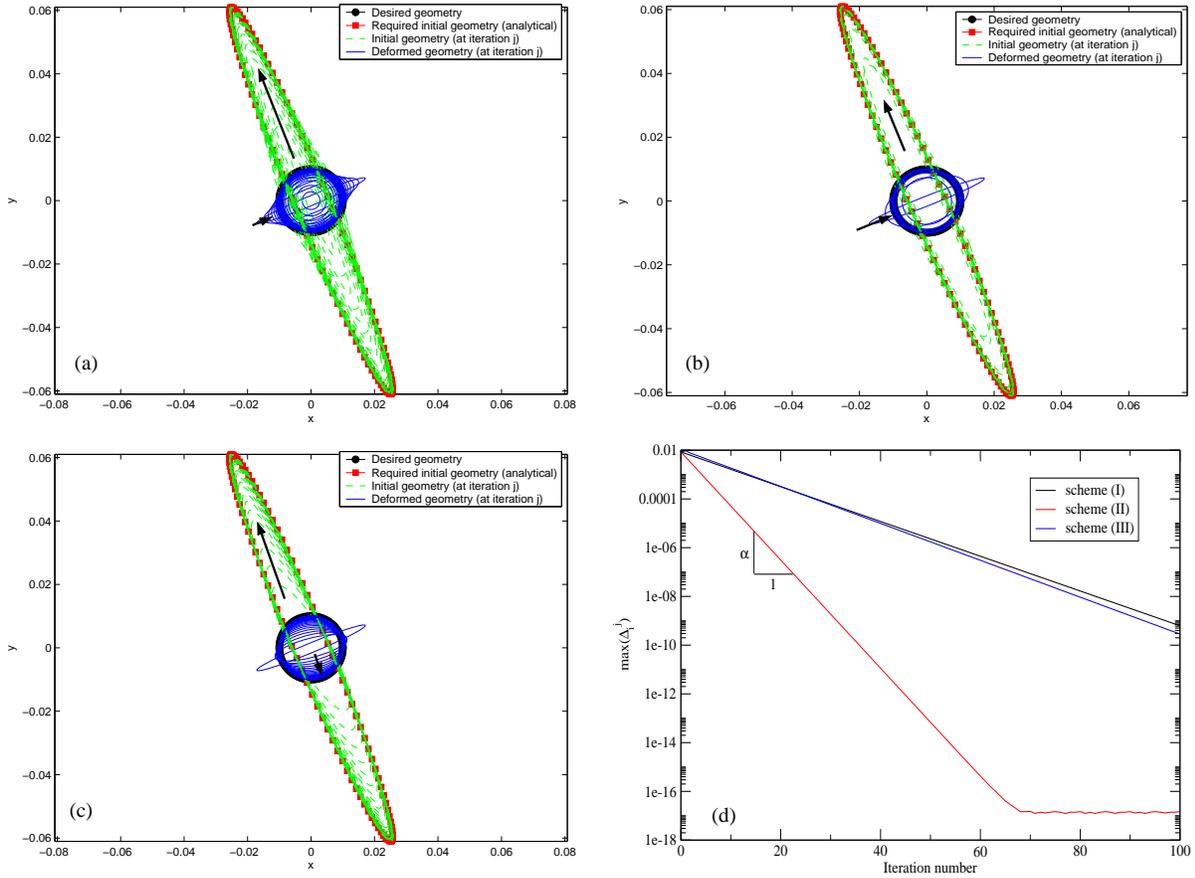}}
\caption{\textbf{(a)}, \textbf{(b)} and \textbf{(c)}: Initial and deformed geometries using the \textbf{schemes (I), (II)} and \textbf{(III)} respectively, $\alpha$=0.6 and 100 boundary points. The arrows indicate the direction of increasing $j$. The figures also show the desired geometry and the true required initial geometry (analytical solution); \textbf{(d)}: corresponding plot of $\max(||\vec{\Delta}_i^j||) \; , \; i\in [1,L]$ against the iteration number.}
\label{fig:alpha06}
\end{center}
\end{figure}

Results are presented on Figure \ref{fig:alpha06} for $\alpha$=0.6, 100 boundary points ($L=100$) and the three different schemes. Clearly, for the three schemes each iteration produces an initial geometry which comes closer and closer to the true required initial geometry calculated analytically using eqs.\ (\ref{eq:truex}) and (\ref{eq:truey}) and consequently the deformed geometry comes closer and closer to a disc of radius 0.01. The naked eye can readily distinguish that fewer iterations are necessary with \textbf{scheme (II)} (see Figure \ref{fig:alpha06} \textbf{(b)}) to identify the required initial geometry than with the two other schemes (Figures \ref{fig:alpha06} \textbf{(a)} and \textbf{(c)}). This is emphasized on Figure \ref{fig:alpha06} \textbf{(d)} where the maximum of the norm of the residual vector from all the boundary points is plotted against the number of iterations. As predicted by the previous analysis, the convergence rate defined by the ratio of the norm of the residual vector at the current iteration to the one at the previous one, i.e. $||\vec{\Delta}_i^j||/||\vec{\Delta}_i^{j-1}||$, is precisely equal to $\alpha$=0.6 for \textbf{scheme (II)}. This is simply a result of eqs.\ (\ref{eq:convscheme}). It is worse for \textbf{scheme (I)} and \textbf{scheme (III)}. The convergence rate is equal to 0.8486 in the former case and 0.84 in the latter.

\begin{figure}
\begin{center}
\scalebox{0.7}{\includegraphics{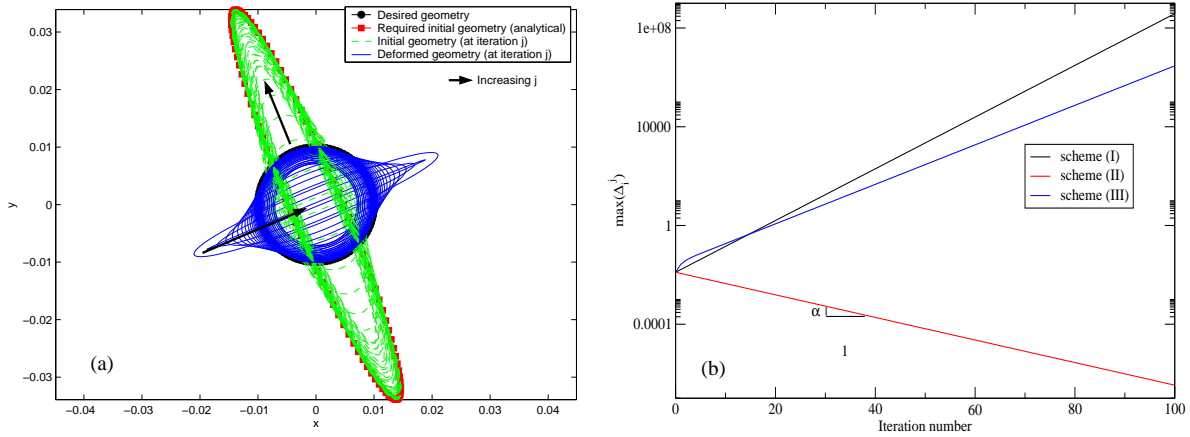}}
\caption{\textbf{(a)}: Initial and deformed geometries using \textbf{scheme (II)}, $\alpha$=0.9 and 100 boundary points. The arrows indicate the direction of increasing $j$. The figure also shows the desired geometry and the true required initial geometry (analytical solution); \textbf{(b)}: corresponding plot of $\max(||\vec{\Delta}_i^j||) \; , \; i\in [1,L]$ against the iteration number.}
\label{fig:alpha09}
\end{center}
\end{figure}

The true benefit of \textbf{scheme (II)} becomes even greater when the value of $\alpha$ is increased. As shown on Figure \ref{fig:alpha09} \textbf{(b)}, \textbf{scheme (I)} and \textbf{scheme (III)} fail to converge completely for $\alpha$=0.9 whereas \textbf{scheme (II)} successfully converges towards the required initial shape. Again, in virtue of the previous analysis, the convergence rate is precisely equal to 0.9. The calculated initial and deformed geometries along with the desired one and the corresponding analytical solution are shown on Figure \ref{fig:alpha09} \textbf{(a)}. 

\begin{figure}
\begin{center}
\scalebox{0.7}{\includegraphics{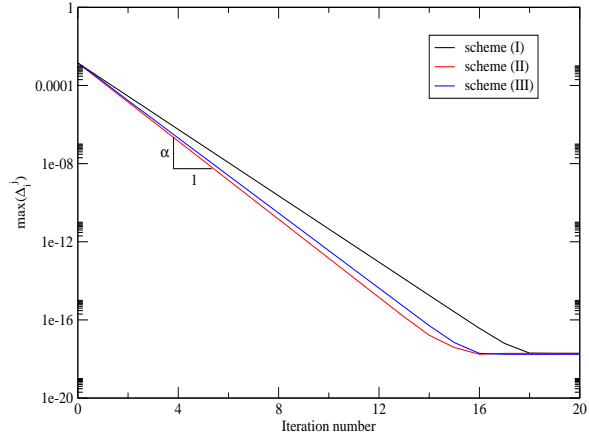}}
\caption{Plot of $\max(||\vec{\Delta}_i^j||) \; , \; i\in [1,L]$ against the iteration number for the three schemes and $\alpha$=0.1.}
\label{fig:alpha01}
\end{center}
\end{figure}

Finally, since \textbf{scheme (III)} is expected, according to the previous analysis, to be well suited for smaller deformations, an additional test is performed with $\alpha =0.1$. The results, plotted on Figure \ref{fig:alpha01}, show a rapid convergence of the three schemes towards the required initial geometry. The convergence rate of \textbf{scheme (III)} is equal to 0.11 and is close to the value of 0.1 for \textbf{scheme (II)}. \textbf{Scheme (I)} displays the worst convergence rate with a value of 0.1414 which confirms the benefit of adding at least an approximate corrective term in the updating of the boundary node location (eq.\ (\ref{eq:sub4})).  
\section{Conclusions}
This paper discusses the possibility to identify the required initial shape so that the deformed one matches very precisely a prescribed shape by means of an iterative scheme. It consists in using the desired work-piece geometry as an initial guess to the required initial geometry and back-tracking the latter by iteratively updating the locations of a set of boundary points. Based on an analysis of the convergence properties of the scheme, a corrective vector for the updating of the boundary point locations is derived. The addition of the latter is shown to extend the applicability of the scheme explored in \cite{Sellier04} to arbitrary deformations when the shear and volumetric strains are equally important. The convergence is proven to be unconditional for small deformations and in the limit of the validity of the Taylor series expansion. Moreover, as illustrated by the test problem, significant improvement in the convergence rate is achieved by introducing this corrective term and to a lesser extent when the approximated corrective term is used.

This scheme could offer a valuable alternative to other approaches for initial shape identification based on sensitivity analysis and optimization methods. At least four potential benefits may be outlined. No parameterization of the work-piece geometry is required since an arbitrarily large number of boundary points can be selected. The method is purely geometric and therefore its success does not depend on the type of constitutive law. The method can easily be used in combination with a commercial code for the computation of the displacement field, see \cite{Sellier04} for an example. The question of how closely the deformed geometry matches the desired one which is not necessarily a trivial one is easily answered thanks to the introduction of the residual vector whose norm gives a clear measure of the mismatch between the geometries. This residual vector can only be defined because the set of reference points (at the desired locations) is defined at the first iteration of the scheme when the desired geometry is chosen as a first guess for the required initial geometry.
\ack The author gratefully acknowledges the funding of the European Union through the MAGICAL project.       


\section*{References} 
\begin{thebibliography}{99}
\bibitem{Chung03} Chung S H, Fourment L, Chenot J-L and Hwang S M 2003 {\it Int. J. Numer. Meth. Engng.} {\bf 57} 1431--44
\bibitem{Vielledent01} Vieilledent D and Fourment L 2001 {\it Int. J. Numer. Meth. Engng.} {\bf 52} 1301--21
\bibitem{Sousa02} Sousa L C, Castro C F, Antonio C A C and Santos A D 2002 {\it J. Mater. Process. Technol.} {\bf 128} 266--273
\bibitem{Sprekels} Sprekels J, Goldberg H and Troeltzsch F {\it Preprint} DFG-Preprint series ``Anwendungsbezogene Optimierung und Steuerung", Report No. 520
\bibitem{Park83} Park J J, Rebelo N and Kobayashi S 1983 {\it Int. J. Mach. Tool D. R.} {\bf 23} 71--79
\bibitem{Kim90} Kim N and Kobayashi S 1990 {\it Int. J. Mach. Tools Manufact.} {\bf 30} 243--68
\bibitem{Zhao94} Zhao G, Wright E and V. Grandhi R {\it Int. J. Mach. Tools Manufact.} {\bf 35} 1225--39
\bibitem{Sellier04} Sellier M 2004 {\it submitted to Int. J. Form. Proc.}
\end{thebibliography}
\end{document}